\documentclass[12pt]{amsart}

\usepackage{amscd}
\usepackage{amsmath,amsfonts,amsthm,epsfig,latexsym,graphicx,amssymb,psfrag}
\usepackage[english]{babel}
\newtheorem{definition}{Definition}[section]
\newtheorem{theorem}{Theorem}
\newtheorem{lemma}{Lemma}[section]

\newtheorem{rema}{Remark}

\newtheorem{corollary}{Corollary}
\newtheorem{proposition}{Proposition}

\newtheorem{rem}{Remark}[section]

\newfont{\df}{cmssbx10}

\def\R{I\kern -0.37 em R}
\def\N{I\kern -0.37 em N}
\def\Z{I\kern -0.37 em Z}

\def\supess_#1{\mathop{\rm supess}\limits_{#1}}
\def\infess_#1{\mathop{\rm infess}\limits_{#1}}

 \def\NN{{\mathbb N}} 

 \def\RR{{\mathbb R}} 
\def\TT{{\mathbb T}}

 \def\ZZ{{\mathbb Z}}

\begin{document}

\title[Burnside problem on $\mathbb{T}^2$.]{\bf Burnside problem for measure preserving groups  {and for 2-groups}  of toral homeomorphisms.}
\author{Nancy Guelman and Isabelle Liousse} \thanks{This paper was partially
supported by Universit\'{e} de Lille 1, PEDECIBA, Universidad de la
Rep\'{u}blica and IFUM}

\address{{\bf  Nancy Guelman}
IMERL, Facultad de Ingenier\'{\i}a, Universidad de la Rep\'ublica,
C.C. 30, Montevideo, Uruguay.  \emph{nguelman@fing.edu.uy}.}

\address{{\bf    Isabelle Liousse}, UMR CNRS 8524, Universit\'{e} de Lille1,
59655 Villeneuve d'Ascq C\'{e}dex,   France.  \emph {liousse@math.univ-lille1.fr}. }

\begin{abstract} A group $G$ is said to be periodic if for {every} $g\in G$ there exists a positive integer $n$  with $g^n=\mathrm{Id}$.  We prove that a finitely generated periodic group of homeomorphisms on the 2-torus
that preserves a probability measure $\mu$ is finite. Moreover if the group consists of homeomorphisms isotopic to
the identity, then it is  abelian and  acts freely on $\mathbb{T}^2$. In the Appendix, we show that {every}  finitely generated 2-group of toral homeomorphisms is finite.

\end{abstract}


\maketitle
\section {Introduction.}

\begin{definition}
A group $G$ is said to be {\bf periodic} if  {every} $g\in G$ has finite
order, that is, there exists a positive integer $n$  with $g^n=\mathrm{Id}$.

\end{definition}

One of the oldest problem in group theory was first posed by William
Burnside in 1902 (see \cite{B}): {\em ``Let $G$ be a finitely generated periodic group.
Is $G$ necessary a finite group?" }

It is obvious that an abelian finitely generated periodic group is
finite.

In 1911, Schur (see \cite{Sh}) proved that this is true for
subgroups of $\mathrm{GL}(k,\mathbb C)$, $k\in \mathbb N$.

\smallskip

But, in general,  according to Golod (see \cite{Go}) the answer is negative.  Adjan and Novikov  (see \cite{AN})  improved this result when  the orders are bounded.

Later,  Ol'shanskii, Ivanov  and  Lysenok  (see \cite{Ol}, \cite{Iv} and \cite{Ly}) exhibited many examples of infinite,  finitely generated and periodic groups  with bounded orders.

\medskip
One of the most interesting examples is the Grigorchuk group, $\Gamma$. It is
a subgroup of the automorphism group of the binary rooted tree $T_2$. It is
generated by   four specific $T_2$-automorphisms $a,b,c,d$, satisfying that
any  of the elements $a,b,c,d$ has order 2 in $\Gamma$, {so that every} element
of $\Gamma$  can be written as a positive word in $a,b,c,d$ without using inverses.
It has been proved (see \cite{Gr}) that the Grigorchuk group is infinite, and it is a 2-group,  that is, every element in $\Gamma$ has finite order that is a power of 2. This shows that the Grigorchuk group is a finitely generated infinite group satisfying that
every element has finite order.

\medskip

The problem raised by Burnside is still open for groups of homeomorphisms (or
diffeomorphisms) on closed manifolds. Very few examples are known.

\medskip

{This question for groups of homeomorphisms and the  following example  are pointed out in \cite{RS} and } was communicated to us by Navas:

A non trivial circle homeomorphism of finite order has no fixed
points, and then a periodic group acting on a circle acts freely.
Moreover, H\"older theorem  states that a group acting freely on the
circle is abelian (see, for example, \cite{RS} or section 2.2.4 of \cite{Na}). If it is also finitely generated then it is
finite.

Therefore, it holds that  {\em  ``a finitely generated periodic  group of
circle homeomorphisms is finite".}

Finally, we note that, even in the circle case,  the hypothesis on finiteness of the generating set is crucial: the group
consisting of all rational circle rotations is periodic and infinite.

\medskip

Rebelo and Silva (see \cite{RS}) proved that any finitely generated
periodic subgroup of $C^2$-symplectomorphisms of a  compact
4-dimensional symplectic manifold is  finite, provided that
the fundamental class in $H^4 (M, \ZZ)$ is a product of classes in $
H^1 (M, \ZZ)$. Another result in \cite{RS} is the following:

Let $ M$  be a compact (oriented) manifold whose fundamental class
in \\ $H^n (M, \ZZ)$  is a product of elements  in  $ H^1 (M, \ZZ)$ and
whose mapping class group  is finite. Then, any finitely generated subgroup G of $\mathrm{Diff }_\mu ^2 (M )$ whose
elements have finite order is finite, where  $\mu$ is a probability
measure on M and $ {\mathrm{Diff}}_\mu ^2(M)$ is the subgroup of
orientation-preserving $C^2$ -diffeomorphisms of M preserving $\mu$.

\bigskip

In this paper, we study a related question : we consider  finitely
generated periodic groups of homeomorphisms of the 2-torus,  $\TT^2$.
 {Our first result is an extension from the $C^2$ to the $C^0$-case of the previous result of \cite{RS}: }

\begin{theorem}
\label{theo}
A finitely generated periodic group of  homeomorphisms
of $\TT^2$ that preserves a probability  measure $\mu$ is finite.

Moreover, if the group consists of homeomorphisms isotopic to
the identity  then { it is conjugate  to a group of rational toral translations (see  Definition \ref{trans})
(in particular it is abelian and acts freely on $\TT^2$)}.

\end{theorem}

\begin{corollary}\label{coramena}
 An amenable finitely generated periodic group of  homeomorphisms
of $ \TT^2$  is finite and if it  consists  {of}  homeomorphisms isotopic to
the identity, {in which case it is also} abelian and acts freely on $\TT^2$. In other words, if an amenable finitely generated periodic group acts on $\TT^2$ by homeomorphisms isotopic to the identity, then the action factors through a finite abelian group which acts freely on $\TT^2$.
\end{corollary}

\begin{rem}
{Since every finite group is amenable, Corollary \ref{coramena} yields a new }  proof of Lemma 4.1 of a recent paper of Franks and Handel (see \cite{FH}). It claims
that "If G is a finite group which acts on $\TT^2$ by homeomorphisms isotopic to
the identity, then the action factors through an abelian group which acts freely on $\TT^2$."
\end{rem}
Since the Grigorchuk group is amenable (see \cite{Gr}) we have the following:
\begin{corollary}\label{grig}
 The Grigorchuk group can not act faithfully on  $ \TT^2$.
\end{corollary}

Finally, in the Appendix, we prove a more general statement of Corollary \ref{grig}. {More precisely,  for 2-groups of toral homeomorphisms, we  are able to dispense with the assumption concerning the existence of an invariant probability measure.}

\begin{theorem}\label{2-group}
 {Every} finitely generated 2-group of toral homeomorphisms is finite.
\end{theorem}

\section {Preliminaries.}

\begin{definition}\label{trans}
Let $v \in {\mathbb R}^2$,  denote by $ \tilde T_v : {\mathbb R}
^2\rightarrow \mathbb R^2$  the {\bf translation} by $v$: $\tilde T_v(x,y)
= (x,y) +v$ and by $ T_v : \TT ^2\rightarrow \TT^2$ its
induced map on $\TT^2$  which is called {\bf toral translation}.

A toral homeomorphism is called  {\bf pseudo-translation} if it is 
{conjugate} to a toral translation.
\end{definition}

\subsection{Homology representation} \

\noindent Let $\displaystyle \mathrm{Homeo}_{\ZZ ^2}(\RR ^2)$ be the set of homeomorphisms $F: \RR
^2 \rightarrow \RR ^2$ such that $F(\ZZ^2 )\subseteq \ZZ ^2$

\smallskip

\noindent and $\mathrm{Homeo}^0_{\ZZ ^2}(\RR ^2)$ be the set of homeomorphisms $F: \RR
^2 \rightarrow \RR ^2$ such that $F(\ZZ^2 )\subseteq \ZZ ^2$ and $F(x + P) = F(x) +P$, for all $x\in \RR^2$
and $P\in \ZZ ^2$.

\smallskip

\noindent Note that   a 2-torus homeomorphism is
isotopic to identity  if and only if any lift  belongs to $\mathrm{ Homeo}^0_{\ZZ ^2}(\RR ^2)$.

\medskip

Let $ f: \TT^2 \rightarrow \TT^2$ be a homeomorphism and let $F: \RR^2
\rightarrow \RR^2$ be a lift of $f$.  We can associate to $F$ a linear map
$A_F $  defined by:
$$  F( p + (m,n)) = F( p )  + A_F (m,n), {\text { for any }}  m,n  {\text
  { integers.} }$$

This map  $A_F$  {depends} neither on the integers $m$ and $n$ nor on the lift
$F$ of $f$.  In fact,  $A_F$ is the morphism  induced by $f$ on the first
homology group of $\TT ^2$.  So we will denote $A_f$ for $A_F$. {The following properties can easily be checked.}

\noindent{\bf Properties.}
\begin{enumerate}
\item {\em  The map} $A:  \mathrm{Homeo}{(\TT^2)} \rightarrow \mathrm{ GL}(2,\ZZ)$ {\em defined by }$A(f) =  A_f$ {\em is a morphism of groups.}

\item{\em  A toral homeomorphism } $f$ {\em  is isotopic to identity if and only if } $A_f= \mathrm{Id}$.

\item{\em {Every} pseudo translation is isotopic to identity.}

\end{enumerate}

\subsection{Measure rotation set.}  \

In \cite{MZ},  rotation sets of torus homeomorphisms are introduced by Misiurewicz
and Ziemian.

\begin{definition}
Let  $f$ be a  2-torus homeomorphism
isotopic to identity. We denote by  $\tilde f$ a lift of $f$ to  $\mathbb R ^2$. We call  {\bf measure rotation set} of $\tilde f$  the subset of $\mathbb R ^2$ defined by
$$\rho_{\mathrm {mes} } ( \widetilde{f}):= \{ \rho_{\mu} ( \widetilde{f})=\int_{\TT^2} (\widetilde{f} -\mathrm{Id}) d\mu,
{\text {  where } }\mu {\text { is an $f$-invariant probability}
}\}.$$

Note that the map $\widetilde{f} -\mathrm{Id} : \TT ^2\rightarrow \RR^2$ is well defined, since $(\widetilde{f} -\mathrm{Id } )( x+(p,q)) = (\widetilde{f} -\mathrm{Id} )( x)$, for all $(p,q)\in \ZZ^2$.

\end{definition}

\noindent  {\bf Properties of  the measure rotation set.}
\begin{proposition}\label{ptyRS}
 Let $f$ be a  2-torus homeomorphism  isotopic to the identity and  $\tilde{f}$ be a lift of $f$ to $\RR^2$. Let $\mu$ be an $f$-invariant probability, let $h$ be {an arbitrary} 2-torus homeomorphism.

\begin{enumerate}

\item $\rho_{\mu} (\tilde f^n ) = n \rho_{\mu}(\tilde f)$,

\item $\rho_{\mu} (\tilde f +(p,q)) = \rho_{\mu}(\tilde f ) +(p,q) $, for any $(p,q)\in \ZZ^2$.

\smallskip

\hskip -0,6 truecm For pseudo-translations we have:

\item $\rho_{\mu} (\tilde T_v) = v$ for {every} $T_v$-invariant measure,

\item $\rho_{\mu} (\widetilde{h\circ T_v \circ h ^{-1}})=  A_h(v)+ (p_1,q_1)$, for some  $(p_1,q_1)\in \ZZ^2$, for any ${h\circ T_v \circ h ^{-1}}$- invariant measure. Hence, the measure rotation set of a pseudo-translation is a single vector.

\end{enumerate}
\end{proposition}

\begin{proof}
The first three items are direct consequences of definitions.
For last item note that exists $(p_1,q_1) \in \ZZ^2$ such that  $\rho_{\mu} (\widetilde{h\circ T_v \circ h
^{-1}})=\rho_{\mu} (\widetilde{h}\circ \widetilde{T_v} \circ \widetilde{h}^{-1}) + (p_1,q_1).$

\smallskip

{\bf Case 1: } $A_h=\mathrm{Id}$. We prove that  $\rho_{\mu} (\widetilde{h}\circ \widetilde{T_v} \circ \widetilde{h}^{-1})=v$.

By definition
$$ \rho_{\mu} (\widetilde{h}\circ \widetilde{T_v} \circ \widetilde{h}^{-1}) =  \int_{\TT^2} (\widetilde{h}\circ \widetilde{T_v} \circ \widetilde{h}^{-1} -\mathrm{Id}) d\mu=  \int_{\TT^2} (\widetilde{h}\circ \widetilde{T_v}- \widetilde { h}) d(h_*\mu ),$$ where $h_*\mu$ is defined by  $h_*\mu(B)= \mu(h(B))$ for {every} measurable set $B \subset \TT^2$.

Then
$$ \rho_{\mu} (\widetilde{h}\circ \widetilde{T_v} \circ \widetilde{h}^{-1}) =\int_{\TT^2}(\widetilde{h}\circ \widetilde{T_v}-\widetilde{T_v}+\widetilde{T_v}- \mathrm{Id} + \mathrm{Id} -  \widetilde{h})  d(h_*\mu )=$$
$$ =\int_{\TT^2}(\widetilde{h}\circ \widetilde{T_v}-\widetilde{T_v})d(h_*\mu )+\int_{\TT^2}(\widetilde{T_v}- \mathrm{Id})d(h_*\mu )+\int_{\TT^2}(\mathrm{Id} -  \widetilde{h})  d(h_*\mu ).$$ These integrals are well defined since $h$ is isotopic to identity.

Since $\mu$ is $h\circ T_v \circ h
^{-1}$-invariant, the measure $h_*\mu $ is $T_v$-invariant. Hence,
$$ \rho_{\mu} (\widetilde{h}\circ \widetilde{T_v} \circ \widetilde{h}^{-1}) =\int_{\TT^2}(\widetilde{h} -\mathrm{Id})d(h_*\mu )+\int_{\TT^2}(\widetilde{T_v}- \mathrm{Id})d(h_*\mu )+\int_{\TT^2}(\mathrm{Id}-  \widetilde{h})  d(h_*\mu )=$$ $$= \int_{\TT^2}(\widetilde{T_v}- \mathrm{Id})d(h_*\mu )=v.$$

\begin{rem} Let $M \in \mathrm{GL}(2,\RR)$, we have that $M\circ \widetilde{T_v} \circ M^{-1}=\widetilde{T}_{Mv}$.
\label{rem21}
\end{rem}

{\bf Case 2: } General case.

Since   $A:  \mathrm{Homeo}{(\TT^2)} \rightarrow  \mathrm{GL}(2,\ZZ)$ is a morphism of groups, we have that $$A_{\widetilde h \circ A_h^{-1}}= A_{\widetilde h} \circ A_{A_h^{-1}}=A_h \circ A_h^{-1} $$ due to linearity of $A^{-1}_h$.
Therefore
$A_{\widetilde h \circ A_h^{-1}}=\mathrm{Id}$, so  $\widetilde {h} \circ A_h^{-1} \in {\mathrm{Homeo}}^0_{\mathbb Z ^2}(\RR ^2).$

As a consequence of Remark \ref{rem21}, we have that $$\widetilde{h}\circ \widetilde{T_v} \circ \widetilde{h}^{-1}=\widetilde{h}\circ A_h^{-1} \circ  A_h \circ \widetilde{T_v} \circ  A_h^{-1} \circ  A_h \circ \widetilde{h}^{-1}=\widetilde{h}\circ A_h^{-1} \circ  \widetilde{T}_{A_h v}  \circ  A_h \circ \widetilde{h}^{-1}$$

By case 1, we have that  $\rho_{\mu} (\widetilde{h}\circ A_h^{-1} \circ  \widetilde{T}_{A_h v}  \circ  A_h \circ \widetilde{h}^{-1})= A_h v$.

Hence, $\rho_{\mu} (\widetilde{h\circ T_v
\circ h ^{-1}})=  A_h(v)+ (p_1,q_1)$, for some  $(p_1,q_1)\in
\ZZ^2$. \end{proof}
\medskip

\begin{rem}
Another proof of this property can be {obtained} by using the same property for the classical rotation set (see for example section 3 of \cite{GL2})  and the relations between different rotation sets.

\end{rem}

\begin{corollary}
If $K_0$ is a subgroup of $\mathrm{Homeo} (\TT ^2)$ consisting of pseudo-translations then the {\bf rotation  map}
$\rho : K_0 \rightarrow \TT^2 = \RR^2/ \ZZ^2$ given  by $\rho (k) = \rho_{\mathrm {mes} } (\tilde k)$ (mod $\ZZ^2$) is well defined.
\end{corollary}

\medskip

\bigskip

From now, we consider finitely generated periodic groups of toral homeomorphisms.

\medskip

\section{Classification of finite abelian  groups of  isotopic to identity toral homeomorphisms.}

\medskip

In this section, we begin  giving  a  classification of finite order isotopic to identity homeomorphisms of the torus
up to conjugacy.

\medskip

\begin{proposition}\label{isotop}
A finite order, isotopic to identity,  toral homeomorphism is
{conjugate} to a rational translation. In other words, it is a rational pseudo-translation.
\end{proposition}

\medskip

This result is  well known:
 {Theorem 2.8 of \cite{E} asserts that ''if $f$ is  a finite order  homeomorphism of a compact orientable surface $M$,  then  there is a riemannian metric of constant curvature on $M$  such that $f$ is a diffeomorphism preserving the riemannian metric.'' }

  {By Killing-Hopf theorem, a torus of constant curvature is  isometric to an {euclidean} torus.}  Moreover, an {euclidean} isometry of $\RR^2$ is $x\mapsto Ax + v$, with $v\in \RR^2$, $A\in \mathcal O (2, \RR)$ the linear group consisting of rotations and reflections in lines.

We claim that  an isotopic to identity {euclidean} toral isometry $h$ is a translation. Indeed,  let $H$ be  a lift of $h$ to $\RR^2$, $H(x)= Ax + v$, and for any  vector $P\in \ZZ^2$ we have that  $H(x+P) = Ax + AP +v=  H(x) + A P$.  Since $h$ is isotopic to identity, it follows that $H(x+P)= H(x)+P$. Therefore $A=\mathrm{Id}$ and $H$ is a translation.

Hence,  a finite order, isotopic to identity,  toral homeomorphism is  
{conjugate} to a translation by some rational vector.

A proof of this statement can also be found in  \cite{GLM}.

\begin{definition}
Let $G$ be a group of toral homeomorphisms, we denoted  by $\mathbf {G_0}$ the subgroup $ \{g\in G : g$ is isotopic to identity $\}$.

\end{definition}

As a corollary of Proposition \ref{isotop} and the remark that any pseudo-translation is isotopic to identity, we get

\begin{corollary}\label{pseudo}

If $G$ is a periodic group, then $G_0=\{g\in G : g$ is a pseudo-translation $\}$ and acts freely on $\TT^2$.

In particular, the subset of $G$ consisting of pseudo-translations is a subgroup.
\end{corollary}

Note that a group of toral pseudo-translations always acts freely, since a non trivial pseudo-translation does not admits fix point.

\medskip
We will generalize Proposition \ref{isotop} for a finite abelian group of isotopic to identity toral homeomorphisms. More precisely, we prove

\begin{proposition}\label{group}
A finite abelian group of isotopic to identity toral homeomorphisms is
{conjugate} to a group of  translations.
\end{proposition}

\begin{proof}  Let $G_0$ be finite abelian group of isotopic to identity toral homeomorphisms. We have proved that  {every} element of $G_0$ is {conjugate} to
a translation. We can suppose, up to conjugacy,  that $G_0 = <T_{v_1}, f_2,...,f_p >$.

\medskip

We will argue by induction {on the number of generators}. If $G^1: = <T_{v_1}>$ it is  {obviously   conjugate} to a group of  translations.

\noindent Let $k\in \{1, ..., p-1\}$, we suppose that  $G^k:= <T_{v_1}, f_2,...,f_k>$ is {conjugate} to a group of  translations and we will prove that  $G^{k+1}=<T_{v_1}, f_2,...,f_{k+1}>$ is also {conjugate} to a group of translations.

\noindent Using the inductive hypothesis, we can suppose, up to conjugacy, that $G^{k}=<T_{v_1},...,T_{v_k}>$ and consequently that
$G^{k+1} = <T_{v_1},...,T_{v_k},f_{k+1}>$.

First, we note that the quotient $M=\TT^2 / G^k$ is a torus. Indeed, it is an orientable, compact surface ($G^k$ is finite and acts freely by orientation preserving homeomorphisms) and the fundamental group of $\TT ^2$ injects into
the fundamental group of $M$ so the only possible case is that $M =\TT^2$.

Let us denote by $\pi : \TT^2 \rightarrow \TT^2/ G^k$ the canonical projection.

Since $f_{k+1}$ commutes with any $T_{v_i}$, the map $\bar f_{k+1} : \TT^2/ G^k \rightarrow \TT^2/ G^k  $ such that $\bar f_{k+1} \circ \pi=  \pi\circ  f_{k+1}$
is well defined and it is of finite order and isotopic to identity. Hence, by Proposition \ref{isotop}, there is  a homeomorphism $h :  \TT^2/ G^k \rightarrow \TT^2/ G^k  $ such that $h\circ \bar f_{k+1} \circ h^{-1}$ is a toral translation.

Let $\tilde h$ be a  lift of $h$ on $\TT^2$ ($\tilde h\circ  T_{v_i} \circ \tilde h ^{-1} = T_{w_i}$, where $\displaystyle w_j= \sum_{i=1} ^k  n_iv_i$,  $n_i\in \ZZ$).

It follows  that  $\tilde h\circ  f_{k+1} \circ \tilde h ^{-1}$ is also a  translation.

Finally, $\tilde h \circ G_{k+1}   \circ \tilde h ^{-1}$   is  a group of  translations. \end{proof}

\medskip

\section{Reduction to  groups of  rational toral pseudo-translations.}

\medskip

\begin{proposition}\label{finiteindex}

Let $G$ be a finitely generated periodic subgroup of  toral homeomorphisms.  Then $G_0$, the subgroup of $G$  consisting  of pseudo-translations,  is of finite index in $G$.

\end{proposition}

\begin{proof} Denote by $A : G \rightarrow \mathrm{GL}(2,\ZZ)$  the representation of $G$ in the homology of $\TT^2$, induced by the map $A$ defined in section 2.1. Note that $G_0= A^{-1}(\mathrm{Id})$.

{Since} $G$ is a finitely generated periodic group, its  image  by $A$, the group $A(G)$, is a finitely generated  periodic subgroup  of $\mathrm{GL}(2,\ZZ)$. Then, by Schur's result (\cite{Sh}),  $A(G)$ is finite.

Hence, since the quotient group $G/A^{-1}(\mathrm{Id})$ is isomorphic to $A(G)$, $G_0= A^{-1}(\mathrm{Id})$ has finite index in $G$. This ends the proof. \end{proof}

As a direct consequence of this proposition we have that:

\begin{corollary} \label{index} $G$ is finite if and only if  $G_0$ is finite.
\end{corollary}

For proving our main theorem, it is enough to prove that $G_0$ is finite, in fact, we will prove that $G_0$ is  a finitely generated  periodic  abelian group.

\medskip

\section{Burnside problem for groups of rational toral  pseudo-translations.}

\medskip

Note that a rational toral  pseudo-translation is of finite order, then groups of rational toral  pseudo-translations are periodic groups.

We note that as in the circle case, such a group acts freely on the torus, but the H\"older theorem does not hold on the torus and the proof on the circle given in the introduction, {cannot} be adapted to the torus.

 We add the hypothesis that the rotation map is a morphism.

\smallskip

\begin{proposition}\label{translation}

Let $K_0$ be a group of  toral  rational pseudo-translations. Suppose that the rotation map defined on $K_0$ is a morphism into $(\TT^2, +)$. Then $K_0$ is  abelian.

Morever, if $K_0$ is  finitely generated  then $K_0$ is finite.
\end{proposition}

\begin{proof}   Since the rotation map is a morphism, the commutator subgroup of
$K_0$ consists of homeomorphisms of trivial rotation vector, indeed
$\rho ([f,g])= \rho (f)+\rho (g)-\rho (f)-\rho (g)=0$. As $K_0$ only
contains pseudo-translations, an element $g_0$ of $[K_0, K_0]$  is
{conjugate} to a translation $T_v$ by a homeomorphism $h$ and $g_0$
has  trivial rotation vector. Since $(0,0)=\rho (g_0)= A_h(v)$ $(mod
\  \ZZ^2)$, then $ A_h(v) \in \ZZ^2$ so $ v\in A_h^{-1}(\ZZ^2) \subset \ZZ^2$. Hence $T_v=\mathrm{Id}$ therefore $g_0$ is {conjugate} to $\mathrm{Id}$,
consequently $g_0$  is trivial. Finally,  $K_0$ is abelian.

Moreover, if $K_0$ is  finitely generated  then,  as it is abelian and periodic, it is  finite. \end{proof}

\section{Proof of the theorem.}

Let $G$ be a finitely generated periodic group preserving  a probability measure $\mu$ and  $G_0$ be the  subgroup consisting of pseudo-translations of $G$.
We begin by proving that the rotation map $\rho$ on $G_0$ is a morphism in order to apply Proposition \ref{translation}. More precisely, we prove a  more general statement:

\begin{proposition}\label{rotation map}
Let $K_0$ be a subgroup of $\mathrm{Homeo} (\TT ^2)$ consisting of pseudo-translations preserving a probability measure $\mu$.  Then the rotation map $\rho : K_0 \rightarrow \TT^2$ is a morphism.

\end{proposition}

\begin{proof} Let $f, g$ in $K_0$. We have to prove that $\rho(f\circ g) = \rho(f) + \rho (g) $.

By definition, $\rho(f\circ g) = \rho_{\mathrm {mes} } ( \widetilde{f\circ g })=\rho_{\mathrm {mes} } ( \widetilde{f} \circ \widetilde{g })$ $(mod \ \ZZ^2$). We have that:

$$\rho_{\mathrm {mes} } ( \widetilde{f}\circ \widetilde{g }):=
 \int_{\TT^2} (\widetilde{f} \circ \widetilde{g} (x) - x ) d\mu(x)= \int_{\TT^2} (\widetilde{f} \circ \widetilde{g}(x)  -\widetilde{g}(x)) d\mu(x)   + \int_{\TT^2} ( \widetilde{g}(x)  -x) d\mu(x)  =$$ $$ = \int_{\TT^2} (\widetilde{f}(y) - y) d(g^{-1}_*\mu )(y)  + \int_{\TT^2} ( \widetilde{g}(x)  -x) d\mu(x),$$  with the change of variables on $\TT^2$,  $y= g(x)$.

\noindent Since $\mu$ is $K_0$-invariant, $g^{-1}_*\mu = \mu $ and therefore $\rho_{\mathrm {mes} } ( \widetilde{f}\circ \widetilde{g })= \rho_{\mathrm {mes} } ( \widetilde{f})+ \rho_{\mathrm {mes} } (\widetilde{g })$.
 \newline
Taking this equality $(mod \ \ZZ^2)$, we get the proposition. \end{proof}

 \bigskip

We can go back to prove main theorem.

According to Proposition \ref{rotation map} with $K_0= G_0$, the rotation map $\rho$ on $G_0$ is a morphism.

On the other hand, by Proposition \ref{finiteindex}, $G_0$ has finite index in $G$. As
$G$ is finitely generated, it follows from Schreier's lemma (see for example \cite{S}), that
states that  {every} subgroup of finite index in a finitely generated
group is finitely generated, that $G_0$  is finitely generated.

Hence, we use Proposition \ref{translation} with $K_0= G_0$, to prove that $G_0$ is abelian and finite.
This  implies that  $G$ is finite.

Moreover, if $G$ consists of homeomorphisms isotopic to the identity, then $G=G_0$ is abelian and consists of pseudo-translations. Also, we have already noted that a group of pseudo-translations acts freely, since a pseudo-translation with a fixed point is necessarily trivial. By Proposition \ref{group} $G$ is {conjugate} to a group of rational translations.

\medskip

\section{Appendix: groups generated by order 2 elements and 2-groups}

Burnside (\cite{B}) noted that a finitely generated periodic group whose elements have order 2 is finite and abelian. This is a consequence of the following

\begin{rema}\label{order}
If $f$, $g$ and $f\circ g$ have order $2$ then $f$ and
$g$ commute.

($\mathrm{Id}= (f\circ g)^2= f\circ g \circ f^{-1}\circ g^{-1}$.)

\end{rema}

\subsection{Groups generated by order 2 elements}

In this section, we prove that a finitely generated periodic group of isotopic to identity toral homeomorphisms,
whose generators have order 2 is finite and isomorphic to $\{1\}$, $\ZZ /2\ZZ$ or $\ZZ /2\ZZ \oplus  \ZZ /2\ZZ$. This is a consequence of the fact that {every} element of order 2 in $G_0$, a periodic group of isotopic to identity toral homeomorphisms,  belongs to the center of $G_0$. It will be proven in Proposition \ref{2-com}.

We begin by recalling the following general facts:

\begin{rema}\label{comm}

 Let $G$ be a group generated by $s$ elements $g_1,..., g_s$ of  finite order,
 an element of $G$ can be written $g= g_1^{p_1}....g_s^{p_s} C$, where $C\in [G,G]$ and $p_i\geq 0$ is bounded by the order of $g_i$.
 So the index of $[G,G]$ in $G$ is finite, moreover it is bounded by the product of the
 orders of $g_1,..., g_s$.

\end{rema}

\begin{rema}\label{rever}
Let $h$ be a finite order isotopic to identity   toral homeomorphism,  if $h$ is {conjugate}    to its inverse, then $h$ has order 2.

\end{rema}
This is a consequence of Proposition \ref{ptyRS}. More precisely, we have that $$ \rho(h)=\rho(h)^{-1} = - \rho(h) \  mod \ \ZZ^2,$$ then
$2 \rho( h)= (0,0) \ mod \ \ZZ^2$. As $ h$ has finite order, it
must have order $2$.
Now we can prove:

\begin{proposition}\label{2-com}  Let  $G_0$  be  a periodic group
of toral homeomorphisms isotopic to identity   and $f\in G_0$ of
 order 2. Then $f$ commutes with {every} element of $G_0$.

\end{proposition}

\begin{proof} Let $g \in G_0$.

{\bf Case 1: $g$ has order two.}

We claim that $f\circ g$ is {conjugate} (by $f$) to its inverse. Indeed,
$(f\circ  g)^{-1} = g^{-1}\circ f^{-1} =g\circ f=
f\circ (f\circ g)\circ f = f\circ (f\circ g)\circ f^{-1} $.

According to  Remark \ref{rever},  $f\circ g$ has order 2. Applying Remark \ref{order}, we {conclude} that $f$ and $g$ commute.

\medskip

{\bf Case 2: general case.} We denote by $n$ the order of $g$.

Let us consider $G=<f,g>$, we will prove that $G$ is abelian.
It is easy to verify that:
$$[G,G] \subset \{g^{p_1}f g^{p_2}f ....f g^{p_s}, p_i\in \NN,  0\leq p_i< n , \sum_i p_i = 0 {\text  { mod }} n \}$$

\vskip -4mm

$$= < g^p f g^{-p} , 0\leq p \leq
n-1>.$$

{Every element}  $g^p f g^{-p}$ for $ 0\leq p < n $ has order 2. Then, according to case 1, {all of them } commute.
Hence  $<g^p f g^{-p} , 0\leq p \leq n-1>$
is a finitely generated periodic abelian group, so  it is
finite. It follows that $[G,G]$ is finite and  using Remark \ref{comm}, we get
that $G$ is finite.

Since  $G$ consists of isotopic to identity toral homeomorphisms and  it is
finite (so it preserves a measure on $\TT^2$),   as a consequence of
Theorem \ref{theo}, $G$ is  also abelian. \end{proof}

\begin{corollary}
Let  $G_0=<f_1,...,f_p, w>$  be  a periodic group
of isotopic to identity toral homeomorphisms, where $f_i$, $1\leq i\leq p$
has order 2. Then $G$ is finite.
\end{corollary}

As a consequence of this corollary and Proposition \ref{group}, we have
\begin{corollary}
Let  $G_0=<f_1,...,f_p>$  be  a periodic group
of isotopic to identity toral homeomorphisms, where $f_i$, $1\leq i\leq p$
has order 2. Then $G_0$ is {conjugate} to one of the following translations groups: $\{\mathrm{Id}\}$, $\{\mathrm{Id}, T_{(\frac{1}{2}, 0)} \}$ or $\{\mathrm{Id}, T_{(\frac{1}{2}, 0)},T_{( 0, \frac{1}{2})},T_{(\frac{1}{2}, \frac{1}{2})} \}$.
\end{corollary}
\medskip

\subsection{2-Groups }

The aim of this section is to prove  that {every} finitely generated 2-group of toral homeomorphisms is finite. In particular,
it gives an alternative proof of Corollary \ref{grig}.

\begin{lemma}\label{comcom}
Let $f$ and $g$ be finite order isotopic to identity toral homeomorphisms.

If $f$ and $g^2$ commute then $[g, f^p]$ and $[f^p, g]$, for $0\leq p\leq order(f)$ have order 2.

\end{lemma}

\begin{proof}
By hypothesis $f^p$ and $g^2$ commute, then
$(f^p g) g = g (g f^p) $, hence $(f^p g f ^{-p} g^{-1}) g f ^p g = g^2f^p$.

Consequently,  $(f^p g f ^{-p} g^{-1}) = g^2f^p g ^{-1}f ^{-p } g^{-1} = g ( gf^p g ^{-1}f ^{-p})  g^{-1}  $.

Finally, $[f^p, g]$ and $[g, f^p]$  are {conjugate} (by $g$).

Since  $[g, f^p]= [f^p, g]^{-1}$, it follows by
Remark \ref{rever} that $[g, f^p]$  and $[ f^p,g]$ have order 2. \end{proof}

Next Proposition will be the main tool to prove Proposition \ref{2-group}.

\begin{proposition}\label{carr}
Let $f$ and $g$ be finite order isotopic to identity toral homeomorphisms.

If $f$ and $g^2$ commute then $f$ and $g$ commute.
\end{proposition}

\begin{proof}  The main part of the proof consists of  proving that the group $<f,g>$ is finite.

Let $h\in <f,g>$, we can write $$h= g^{2k_1+\epsilon_1} f ^{\alpha_1}g^{2k_2+\epsilon_2} f ^{\alpha_2}...g^{2k_s+\epsilon_s} f ^{\alpha_s} g^{2k_{s+1}+\epsilon_{s+1}},$$ with $k_i\in \NN$, $\epsilon _i \in \{0,1\}$ and $\alpha_i\in \NN$, $1\leq \alpha_i < order(f)$.

\medskip

As $g ^{2}$  and therefore $g^{2k_i}$ commute with $f $, we have, up to  change $\alpha_i$ and $s$ if were necessary, that
$$ h= g ^{m} g f ^{\alpha_1} g f ^{\alpha_2}...g  f ^{\alpha_s} g^{\epsilon},$$ with $\epsilon \in  \{ 0,1\}$, $\alpha_i,m\in \NN$, $1\leq \alpha_i < order(f)$ and
 $0\leq m < order(g)$.

\medskip

Since $g=g ^2 g ^{-1} $ and $g^2$ commute with $f$, we can write $h$ as
$$h= g ^n g f ^{\alpha_1} g^{-1}  f ^{\alpha_2} g ...   f ^{\alpha_s} g^{\epsilon},$$  with $\epsilon \in  \{0,1\} $ and $n,\alpha_i\in \NN$, $0\leq n < order(g)$, $1\leq \alpha_i < order(f)$. Then

\medskip

$$h= g ^n \left(g f ^{\alpha_1} g^{-1}  f ^{-\alpha_1} \right)   f ^{\alpha_2+\alpha_1} g ...   f ^{\alpha_s} g^{\epsilon}$$

$$= g ^n \left(g f ^{\alpha_1} g^{-1}  f ^{-\alpha_1} \right)  ( f ^{\alpha_2+\alpha_1} g f ^{-\alpha_2-\alpha_1} g ^{-1} ) g f ^{\alpha_3+\alpha_2+\alpha_1}  ...   f ^{\alpha_s} g^{\epsilon},$$

$$= g ^n \left(g f ^{\alpha_1} g^{-1}  f ^{-\alpha_1} \right)  ( f ^{\alpha_2+\alpha_1} g f ^{-\alpha_2-\alpha_1} g ^{-1} ) g f ^{\alpha_3+\alpha_2+\alpha_1}  ...    f ^{\alpha_s+...+ \alpha_1} g^{\epsilon},$$

with $\epsilon \in  \{0,1\} $ and $n,\alpha_i\in \NN$, $0\leq n < order(g)$, $1\leq \alpha_i < order(f)$.  Hence:

$$h= g ^n [g, f ^{\beta_1}] \  [f^{\beta_2}, g] \  [g, f^{\beta_3}] \  ...   \ [\ \  ,\ \  ] \ f ^{m}  g^{\epsilon},$$

with $\epsilon \in  \{0,1\} $ and $m,n, \beta_i\in \NN$, $1\leq m, \beta_i < order(f)$, $0\leq n < order(g)$.

\medskip

Due to the existence of a finite number of  $[g, f^p]$ and $[f^q, g]$, according to
 Lemma \ref{comcom}   and  Remark \ref{2-com}, $[g, f^p]$ and $[f^q, g]$ have order 2 and commute, then
there is only a finite number of possible $h\in <f,g>$, in other words $<f,g>$ is a finite group.

Using  Theorem \ref{theo}, we get that $<f,g>$ is a also an abelian group, proving the claim of this proposition.\end{proof}

\medskip

\begin{proof} of Proposition \ref{2-group}.

By Corollary \ref{index}, it is enough to prove that a  finitely generated 2-group $G_0$ of isotopic to identity toral homeomorphisms is finite.

Let $f$, $g$ in $G_0$. By hypothesis, $g$ has order $ 2 ^{p+1}$, then $g^{2^p}$ has order $2$, so according to Lemma \ref{2-com}, $g^{2^p}$  commute with $f$.

Applying Proposition \ref{carr}, we obtain that $g^{2^{p-1}}$  commute with $f$. Iterating this argument, we get that $f$ and $g$ commute.

Hence, $G_0$ is an abelian finitely generated periodic group, so it is finite. \end{proof}

\bigskip

\noindent {\bf Acknowledgements.} We are grateful to Andr\'es Navas for proposing to us this subject and for fruitful discussions. We thank Christian Bonatti for suggesting us to add Proposition \ref{group}.  { We thank  the referee for giving relevant and helpful comments, suggestions and corrections.}

\bigskip

\bigskip

 \end{document}